\documentclass[reqno,12pt]{amsart}
\pdfoutput=1
\makeatletter
\let\origsection=\section \def\section{\@ifstar{\origsection*}{\mysection}} 
\def\mysection{\@startsection{section}{1}\z@{.7\linespacing\@plus\linespacing}{.5\linespacing}{\normalfont\scshape\centering\S}}
\makeatother        
\usepackage{amsmath,amssymb,amsthm,mathrsfs}
\usepackage{mathabx}
\usepackage{dsfont} 

\usepackage{xcolor}
\usepackage[backref]{hyperref}
\hypersetup{
    colorlinks,
    linkcolor={red!60!black},
    citecolor={green!60!black},
    urlcolor={blue!60!black}
}

\usepackage[open,openlevel=2,atend]{bookmark}
\usepackage[abbrev,msc-links,backrefs]{amsrefs}
\usepackage{doi}

\renewcommand{\PrintDOI}[1]{\doi{#1}}

\usepackage[T1]{fontenc}
\usepackage{lmodern}
\usepackage{microtype}

\usepackage[english]{babel}

\usepackage{geometry}
\usepackage{enumitem}
\def\rmlabel{\upshape({\itshape \roman*\,})}

\def\alabel{\upshape({\itshape \alph*\,})}

\let\polishlcross=\l
\def\l{\ifmmode\ell\else\polishlcross\fi}

\let\emptyset=\varnothing
\let\setminus=\smallsetminus

\makeatletter
\def\moverlay{\mathpalette\mov@rlay}
\def\mov@rlay#1#2{\leavevmode\vtop{%
   \baselineskip\z@skip \lineskiplimit-\maxdimen
   \ialign{\hfil$\m@th#1##$\hfil\cr#2\crcr}}}
\newcommand{\charfusion}[3][\mathord]{
    #1{\ifx#1\mathop\vphantom{#2}\fi
        \mathpalette\mov@rlay{#2\cr#3}
      }
    \ifx#1\mathop\expandafter\displaylimits\fi}
\makeatother

\let\eps=\varepsilon
\let\epsilon=\varepsilon
\let\theta=\vartheta
\let\rho=\varrho
\let\phi=\varphi

\newtheorem{theorem}{Theorem}
\newtheorem{thm}[theorem]{Theorem} 
\newtheorem{lem}[theorem]{Lemma}   

\newtheorem{conj}[theorem]{Conjecture}

\newtheorem{dfn}[theorem]{Definition} 

\def\cC{\mathcal C}
\def\cD{\mathcal D}
\def\ccD{\mathscr D}

\def\cF{\mathcal F}
\def\cG{\mathcal G}

\def\cS{\mathcal S}
\def\cT{\mathcal T}

\def\expec{\mathds E}
\def\nat{\mathds N}
\def\prob{\mathds P}
\def\real{\mathds R}

\def\hh{\overline{h}}

\def\bk{t}
\def\bj{j}

\def\ci{\{\cC_i\}_{i\in[\bk]}}
\def\ui{\{U_i\}_{i\in[\bk]}}

\def\hi{\{h_i\}_{i\in[\bk]}}
\def\hhi{\{\hh_i\}_{i\in[\bk]}}
\def\ffi{\{f_i\}_{i\in[\bk]}}
\def\Fi{\{F_i\}_{i\in[\bk]}}

\def\dmn{\ccD^{m,N}}
\def\dlm{\ccD^{\ell,m}}
\def\dsl{\ccD^{S,\ell}}

\def\ddmn{\cD^{m,N}}
\def\ddlm{\cD^{\ell,m}}
\def\ddslm{\cD^{S,\ell,m}}

\def\bf{\boldsymbol{f}}
\def\bg{\boldsymbol{g}}
\def\bh{\boldsymbol{h}}
\def\bhh{\boldsymbol{\hh}}

\def\dpf{d_{\bf}^{\partial}}

\def\dph{d_{\bh}^{\partial}}

\def\curr{\text{current}}

\def\seq{$(\cG,n,\Delta)$-sequence }
\def\seqs{$(\cG,n,\Delta)$-sequences }

\begin{document}

\title{Packing minor-closed families of graphs
into complete graphs}

\author[Silvia Messuti]{Silvia Messuti}
\address{Fachbereich Mathematik, Universit\"at Hamburg, Hamburg, Germany}
\email{\{\,silvia.messuti\,|\,schacht\,\}@math.uni-hamburg.de}

\author[Vojt\v ech R\"odl]{Vojt\v ech R\"odl}
\address{Department of Mathematics and Computer Science, Emory
  University, Atlanta, USA}
\email{rodl@mathcs.emory.edu}

\author[Mathias Schacht]{Mathias Schacht}

\thanks{The first author was supported through a Doctoral completion scholarship of the University of Hamburg, the second author was supported by NSF grant DMS~1301698, and the third author
was supported through the Heisenberg-Programme of the
Deutsche Forschungsgemeinschaft (DFG Grant SCHA 1263/4-2).}

\keywords{Packings, minor-closed families, trees}
\subjclass[2010]{05C70 (primary), 05C51 (secondary)}

\begin{abstract}
Motivated by a conjecture of Gy\'arf\'as, recently 
B\"ottcher, Hladk\'y, Piguet, and Taraz showed that 
every collection $T_1,\dots,T_t$ 
of trees  on~$n$ vertices with $\sum_{i=1}^te(T_i)\leq \binom{n}{2}$
and with bounded maximum degree,
can be packed into the complete graph on
$(1+o(1))n$ vertices.
We generalise this result where we relax the restriction of packing families of trees to families of graphs of any given non-trivial minor-closed class of graphs.
\end{abstract}

\maketitle

\section{Introduction}

Given graphs~$H$ and~$F$, an \emph{$F$-packing} of $H$ is a collection of edge-disjoint subgraphs of~$H$ that are isomorphic to~$F$.
This definition naturally extends to sequences of graphs. In particular, we say that $\cF=(F_1,\dots,F_t)$ \emph{packs} into $H$ if there exist edge-disjoint subgraphs $H_1,\dots,H_t\subseteq H$ with $H_i$ isomorphic to $F_i$ for every $i\in[t]$.
Gy\'arf\'as' tree packing conjecture~\cite{GL76} initiated a lot of research and asserts the following for the case where $H$ is a complete graph and $\cF$ is a sequence of trees.
\begin{conj}
Any sequence of trees $(T_1,\dots,T_n)$  with $v(T_i)=i$ for $i\in[n]$ packs into~$K_n$.
\end{conj}
The difficulty of this conjecture lies in the fact that it asks for a perfect packing, i.e., a packing where all the edges of $K_n$ are used, since each tree has $e(T_i)=i-1$ edges and hence $\sum_{i\in[n]}e(T_i)=\binom{n}{2}$.
Although some special cases were proven (see, e.g.,~\cite{Ho81} and the references in~\cite{BHPT}), this conjecture is still widely open.

Recently, B\"ottcher, Hladk\'y, Piguet, and Taraz \cite{BHPT} showed that a restricted approximate version holds.
More precisely, they considered a host graph with slightly more than $n$ vertices and trees with bounded maximum degree, while relaxing the assumption on the number of vertices of each tree.

\begin{theorem}[B\"ottcher, Hladk\'y, Piguet, and Taraz]\label{bhpt}
For any $\eps>0$ and any $\Delta \in \nat$ there exists $n_0\in\nat$ such that for any $n\geq n_0$ the following holds for  every $t\in\nat$. 
If $\cT=(T_1,\dots,T_{\bk})$ is a sequence of trees satisfying
\begin{enumerate}[label=\alabel]
\item{$\Delta(T_i)\leq\Delta$ and $v(T_i)\leq n$ for every $i\in [\bk]$, and}
\item{$\sum_{i=1}^{\bk} e(T_i)\leq\binom{n}{2}$,}
\end{enumerate}
then $\cT$ packs into $K_{(1+\eps)n}$.
\end{theorem}
In case $(1+\eps)n$ is not an integer, we should talk about $K_{\lfloor(1+\eps)n\rfloor}$.
However, since we provide asymptotical results, we will omit floors and ceilings here.
The proof of Theorem~\ref{bhpt} is based on a randomized embedding strategy, which draws some similarities to the semirandom nibble method (see e.g.\ \cite{AS08}).
Inspired by the result in \cite{BHPT}, we obtained a somewhat simpler proof of Theorem \ref{bhpt}, which extends from sequences of trees to sequences of graphs contained in any non-trivial minor-closed class.

\begin{thm}\label{mcg}
For any $\eps>0$, $\Delta\in\nat$, and any non-trivial minor-closed family $\cG$ there exists~$n_0\in\nat$ such that for every $n\geq n_0$ the following holds
for every integer $t\in\nat$. 
If~$\cF=(F_1,\dots,F_{\bk})$ is a sequence of graphs from $\cG$ satisfying
\begin{enumerate}[label=\alabel]
\item\label{it:main1}{$\Delta(F_i)\leq\Delta$ and $v(F_i)\leq n$ for every $i\in[\bk]$, and}
\item\label{it:main2}{$\sum_{i=1}^{\bk} e(F_i) \leq\binom{n}{2}$,}
\end{enumerate}
then $\cF$ packs into $K_{(1+\eps)n}$.
\end{thm}
In the following we will consider graphs that do not contain isolated vertices.
In fact, such vertices can easily be embedded after larger components just by picking any vertex of $K_{(1+\eps)n}$ that has not been used before for the same graph.
In the proof we split the graphs~$F_i$ into smaller pieces by removing a \emph{small} separator, i.e., a small subset of the vertex set. We discuss these concepts and a generalisation of Theorem~\ref{mcg} in the next section. 

\section{Main technical result}
We shall establish a generalisation of Theorem~\ref{mcg} for graphs with small separators (see Theorem~\ref{main} below).
In fact, the Separator Theorem of Alon, Seymour, and Thomas \cite{AST90} will provide the connection between Theorem~\ref{mcg} and slightly more general Theorem~\ref{main}.
\begin{thm}[Alon, Seymour, and Thomas]\label{sep}
For every non-trivial minor-closed family of graphs $\cG$ there exists $c_\cG>0$ such that for every graph $G\in\cG$ there exists $U\subseteq V(G)$ with~$|U|\leq c_\cG\sqrt{n}$ such that every component of $G-U$ has order at most $n/2$.
\end{thm}

The graphs we consider in our main result satisfy the following property.
\begin{dfn}\label{deltas}
Given $\delta>0$ and $s\in\nat$, a $(\delta,s)$-separation of a graph $G=(V,E)$ with minimum degree $\delta(G)\geq 1$
is a pair $(U,\cC)$ satisfying 
\begin{enumerate}[label=\rmlabel]
\item $U\subseteq V$, $|U|\leq\delta v(G)$ and
\item \label{sep2} $\cC=G[V\setminus U]$, i.e., the subgraph of $G$ induced on $V\setminus U$, has the property that each component of $\cC$ has order at least two and at most $s$.
\end{enumerate}
We refer to $U$ as the \emph{separator}, and to $\cC$ as the \emph{component graph} of $G$. 
\end{dfn}

Note that, for technical reasons that will become clear later (see equation \eqref{bdfinal}), in~\ref{sep2}. we only allow components of size at least two.
Although the removal of a separator could induce components of size one,
such a separator $U^0$ of $G$ may yield at most $\Delta|U^0|$ components of size one, because in our setting we only deal with graphs $G$ of bounded degree $\Delta(G)\leq \Delta$.
This allows us to add these ``few'' vertices to~$U^0$ without enlarging it too much, and ensure that the resulting set $U$ complies with the definition above.

\begin{dfn}
A family $\cG$ of graphs with minimum degree at least one is $(\delta,s)$-separable if every  $G \in \cG$ admits a
$(\delta,s)$-separation.
\end{dfn}

We will deduce Theorem \ref{mcg} from the following result, in which the condition of~$\cG$ being minor-closed is replaced by the more general property of being $(\delta,s)$-separable.

\begin{thm}\label{main}
For any $\eps>0$ and $\Delta\in\nat$ there exists $\delta>0$ such that for every $s\in\nat$ and any $(\delta,s)$-separable family $\cG$ there exists $n_0\in\nat$ such that for every $n\geq n_0$ the following holds. 
If $\cF=(F_1,\dots,F_{\bk})$ is a sequence of graphs from $\cG$ satisfying
\begin{enumerate}[label=\alabel]
\item{$\Delta(F_i)\leq\Delta$ and $v(F_i)\leq n$ for every $i\in[\bk]$, and}
\item{$\sum_{i=1}^{\bk} e(F_i) \leq\binom{n}{2}$,}
\end{enumerate}
then $\cF$ packs into $K_{(1+\eps)n}$.
\end{thm}

As mentioned above, Theorem \ref{mcg} easily follows from Theorem \ref{main}.
First we show that for any non-trivial minor-closed family $\cG$ and any $\delta>0$ there is some $s$ such that $\cG$ is $(\delta,s)$-separable.
Then we use this fact to deduce Theorem \ref{mcg}.

For a given a graph $G\in\cG$ of order $n$ with minimum degree $\delta(G)\geq 1$ and maximum degree $\Delta(G)\leq \Delta$, we apply Theorem \ref{sep} to all components of $G$ that have some size $r_0$ with ${\frac{n}{2}\leq r_0 \leq n}$.
Since there are at most two such components, at most two applications of Theorem \ref{sep} lead to a separator of size at most $2c_\cG n^{1/2}$ and a set of components all of which have order less than $n/2$.
We then apply Theorem \ref{sep} to all components of $G$ that have some size~$r_1$ with $\frac{n}{4}\leq r_1 < \frac{n}{2}$ and obtain another separator of size at most $4c_\cG \left(\frac{n}{2}\right)^{1/2}$.
At this point all components have order less than $n/4$.
Again, we apply Theorem \ref{sep} to all components of some size $r_2$ with $\frac{n}{8}\leq r_2 < \frac{n}{4}$, add at most $8c_\cG \left(\frac{n}{4}\right)^{1/2}$ more vertices to the separator, and so on.
After $i>0$ such iterations we obtain a separator $U^0\subseteq V(G)$ such that 
\[
|U^0|\leq 2c_\cG n^{1/2}+4c_\cG \left(\frac{n}{2}\right)^{1/2}+\dots+2^{i}c_\cG \left(\frac{n}{2^{i-1}}\right)^{1/2}<
2c_\cG n^{1/2}\cdot\frac{\sqrt{2}^{\,i}-1}{\sqrt{2}-1}<6c_\cG n^{1/2}\,2^{i/2}
\]
and each component of $G-U^0$ has order at most $n/2^i$.
For given $\delta>0$ we can apply this with
\[
i=\left\lfloor2 \log_2 {\left(\frac{\delta n^{1/2}}{6c_\cG(\Delta+1)}\right)}\right\rfloor
\]
and obtain a separator $U_0$ of order at most $\delta n/(\Delta+1)$, and a set of components all of which have order at most $72c_\cG^2(\Delta+1)^2/\delta^2$. 
Note that some of the components in $G-U^0$ may have size one.
However, owing to the maximum degree of $G$ there are at most $\Delta|U^0|$ such components.
By defining $U$ as the separator of size at most $\delta n$ obtained from $U^0$ by adding all these degenerate components of order one, we have shown that the non-trivial minor-closed family~$\cG$ is $(\delta,s)$-separable with $s=72c_\cG^2(\Delta+1)^2/\delta^2$.
Applying Theorem \ref{main} with this~$s$ yields Theorem~\ref{mcg}.

The rest of this paper is devoted to the proof of Theorem \ref{main}.
In Section \ref{prfth} we introduce some definitions and state two technical lemmas that are used in the proof of the theorem, which is given at the end of the section.
Resolvable and almost resolvable decompositions, which we will use to construct our packing, are introduced in Section \ref{dec}.
Finally, the two technical lemmas, Lemma \ref{l1} and Lemma \ref{l2}, are proved in Sections \ref{L1} and \ref{L2}, respectively.
\section{Proof of the main result}\label{prfth}
The following notation will be convenient.
\begin{dfn}
Let $\cG$ be a family of graphs. A $\bk$-tuple of graphs $\cF=(F_1,\dots,F_{\bk})$ with $F_i\in\cG$ and $i\in[\bk]$ is called a \seq if 
\begin{enumerate}[label=\alabel]
\item{$\Delta(F_i)\leq\Delta$ and $v(F_i)\leq n$ for every $i\in[\bk]$, and}
\item\label{it:seqb}{$\sum_{i=1}^{\bk} e(F_i) \leq\binom{n}{2}$.}
\end{enumerate}
\end{dfn}
We will consider \seqs with the following additional properties:
\begin{itemize}
\item $\cG$ will be a $(\delta,s)$-separable family and
\item each graph $F_i$ will be associated with a fixed $(\delta,s)$-separation $(U_i,\cC_i)$.
\end{itemize}
Note that, since we are only considering graphs $F_i$
that do not contain isolated vertices, we have $v(F_i)\leq 2e(F_i)$ and, hence,
\[
\sum_{i=1}^{\bk} |U_i|\leq\sum_{i=1}^{\bk} \delta v(F_i)\leq \delta \sum_{i=1}^{\bk} 2e(F_i) \overset{\text{\ref{it:seqb}}}{\leq} 2\delta\binom{n}{2}< \delta n^2\,.
\]
For a simpler notation we will often suppress the dependence on $U_i$ when we refer to a \seq $(F_1,\dots,F_{\bk})$, since the separator $U_i$ will be always clear from the context.
In a component $C$ from $\cC_i$ we distinguish the set of vertices that are connected to the separator $U_i$ and refer to this set as the \emph{boundary} $\partial C$ of $C$
\[
\partial C=V(C)\cap N(U_i),
\]
where as usual $N(U_i)$ denotes the union of the neighbours in $F_i$ of the vertices in $U_i$.

Moreover, for a component graph $\cC_i$ we consider the union of the boundary sets of its components and set
\[
\partial \cC_i=\bigcup\{\partial C \colon  C \text{ component in } \cC_i\}.
\]
Note that 
\begin{equation}\label{bndr}
|\partial\cC_i|\leq\sum_{u\in U_i} d(u)\leq|U_i|\Delta\leq\delta n \Delta\,.
\end{equation}

For the proof of Theorem \ref{main} we shall pack a given \seq into $K_{(1+\eps)n}$.
The vertices of the host graph $K_{(1+\eps)n}$ will be split into a large part $X$ of order $(1+\xi)n$ for some carefully chosen $\xi=\xi(\eps,\Delta)>0$, and a small part $Y=V\setminus X$.
We will pack the graphs $\ci$ into the clique spanned on $X$ and the sets $\ui$ into $Y$.
For this, we shall ensure that the vertices representing the boundary $\partial\cC_i$ will be appropriately connected to the vertices representing the separator $U_i$.
Having this in mind we will make sure that each vertex of $X$ will only host a few boundary vertices.
In fact, since every edge of the complete bipartite graph induced by $X$ and $Y$ can be used only once in the packing, each vertex $x\in X$ can be used at most $|Y|$ times as boundary vertex for the packing of the sequence $\ci$.
This leads to the following definition.
\begin{dfn}\label{xinbp}
For every $i\in[\bk]$, let $F_i=(V_i,E_i)$ be graphs with separators $U_i,\subseteq V_i$ and component graphs $\cC_i=F_i[V_i\setminus U_i]$.
For a family of injective maps $\bf=\ffi$ with~$f_i\colon V(\cC_i)\rightarrow X$ and for $x\in X$ we define the \emph{boundary degree} of $x$ with respect to $\bf$ by 
\[
\dpf(x)=|\{i\in[\bk]\colon f_i^{-1}(x)\in\partial \cC_i\}|\,.
\]
We call such a family of maps \emph{$b$-balanced} for some $b\in\real$ if  $\dpf(x)\leq b$ for every $x\in X$.
\end{dfn}

Theorem \ref{main} follows from Lemma \ref{l1} and Lemma \ref{l2} below.
Lemma \ref{l1} yields a balanced packing of the component graphs $\ci$ into the clique spanned by $X$ with $|X|\leq(1+\xi)n$.

\begin{lem}\label{l1}
For any $\xi>0$ and $\Delta\in\nat$ there exists $\delta>0$ such that for every $s\in\nat$ and any $(\delta,s)$-separable family $\cG$ there exists $n_0\in\nat$ such that if $\cF$ is a \seq with~$n\geq n_0$, then there exists a $(\xi n)$-balanced packing of the component graphs $\ci$ of all members of $\cF$ into $K_{(1+\xi)n}$. 
\end{lem}

Once we have a balanced packing of $\ci$ into $K_{(1+\xi)n}$, the next lemma allows us to extend it to a packing of $\cF=(F_1,\dots,F_{\bk})$ into a slightly larger clique of size $(1+\eps)n$.

\begin{lem}\label{l2}
For any $\eps>0$ and $\Delta\in\nat$, there exist $\xi>0$ and $\delta>0$ such that for every $s$ and any $(\delta,s)$-separable family $\cG$ there exists $n_0$ such that for any $n\geq n_0$ the following holds.
Suppose there exists a $(\xi n)$-balanced packing of the component graphs~$\ci$ associated with a \seq $\cF$ into $K_{(1+\xi)n}$.
Then there exists a packing of $\cF$ into~$K_{(1+\eps)n}$.
\end{lem}

We postpone the proofs of Lemma \ref{l1} and Lemma \ref{l2} to Section \ref{L1} and Section \ref{L2}. 
Here we describe the proof of our main Theorem based on these two lemmas.

\subsection{Proof of Theorem \ref{main}}
We will first fix all involved constants.
Note that Theorem \ref{main} and Lemma \ref{l2} have a similar quantification.
Hence, for the proof of Theorem \ref{main}, we may apply Lemma \ref{l2} with $\eps$ and $\Delta$ from Theorem \ref{main} and obtain $\xi$ and $\delta'$.
Then Lemma \ref{l1} applied with $\xi$ and $\Delta$ yields a constant $\delta''$.
For Theorem \ref{main} we set $\delta=\min\{\delta',\delta''\}$.
After displaying~$\delta$ for Theorem \ref{main} we are given some $s\in\nat$ and a $(\delta,s)$-separable family $\cG$.

With constants chosen as above, we can apply Lemma \ref{l1} for a \seq $\cF$ which then asserts that the assumptions of Lemma \ref{l2} are fulfilled.
Finally, the conclusion of Lemma \ref{l2} yields Theorem \ref{main}.
\qed

\section{Resolvable and almost resolvable decompositions}\label{dec}

The idea of the proof is to split the  given graphs $(F_1,\dots,F_{\bk})$ into small components, group such components by isomorphism types, and pack components from the same group into complete subgraphs of $K_{(1+\eps)n}$. 
For that we will use Theorem \ref{RCW} and Theorem \ref{R}.

A $K_m$-factor of $K_n$ is a collection of $\frac{n}{m}$ vertex disjoint cliques of order $m$, and a resolvable $K_m$-decomposition of $K_n$ is a collection of
\[
\frac{\binom{n}{2}}{\binom{m}{2}}\frac{m}{n}=\frac{n-1}{m-1}
\]
edge disjoint $K_m$-factors.
Theorem \ref{RCW} states that the obvious necessary divisibility conditions for the existence of a $K_m$-decomposition of $K_n$ are actually sufficient.
\begin{theorem}[Ray-Chaudhury and Wilson]\label{RCW}
For every~$m\geq 2$ there exists~$n_0$ such that if~$n\geq n_0$ and~$n\equiv m\pmod{m(m-1)}$, then~$K_n$ admits a resolvable~$K_m$-decomposition.
\end{theorem}

For general $F$, resolvable decompositions do not necessarily exists (for example it is easy to see that there is no~$n$ for which resolvable~$K_{1,3}$-decompositions of~$K_n$ exist).
Therefore, instead of $F$-factors, we consider $F$-matchings, i.e., sets of vertex disjoint copies of $F$.
\newpage

\begin{dfn}\label{fctr}
A \emph{$(F,\eta)$-factorization} of $K_\ell$ is a collection of $F$-matchings of $K_\ell$ such that 
\begin{enumerate}[label=\rmlabel]
\item each matching has size at least $(1-\eta)\frac{\ell}{v(F)}$, and \label{R1}
\item these matchings together cover all but at most $\eta\binom{\ell}{2}$ edges of $K_\ell$. \label{R2}
\end{enumerate}
\end{dfn}
From these two properties we deduce that the number $t$ of $F$-matchings in an $(F,\eta)$-factorization satisfies 
\[
(1-\eta)\frac{(\ell-1)v(F)}{2e(F)}\leq t\leq\frac{(\ell-1)v(F)}{2e(F)}.
\]
Also note that any $(F,0)$-factorization of $K_\ell$ is a resolvable $F$-decomposition of $K_\ell$.
We will then use the following approximate result, which can be deduced from \cite{FR85} and \cite{PS89} (see also~\cite{AY00}).
\begin{theorem}\label{R}
For every graph $F$ and $\eta>0$ there exists $\ell_0$ such that for every $\ell\geq \ell_0$ there exists an $(F,\eta)$-factorization of $K_\ell$.
\end{theorem}

\section{Packing the components}\label{L1}

The crucial part in the proof of Theorem \ref{main} is Lemma \ref{l1}, which we are going to prove in this section.
In Lemma \ref{l1} we are given a \seq $(F_1,\dots,F_{\bk})$ of graphs from a~$(\delta, s)$-separable family $\cG$ with fixed separations $(U_i,\cC_i)$ associated with each $F_i$.
Our goal will be to construct a $(\xi n)$-balanced packing of the component graphs $\ci$ into~$K_N$, with~$N=(1+\xi)n$. 

The packing of $\ci$ will make use of a resolvable $K_m$-decomposition of $K_N$ (actually we will use a somewhat more complicated auxiliary structure which we will describe in Section~\ref{outass})
and will be realized in  two steps: the \emph{assignment} phase and the \emph{balancing} phase.
\begin{itemize}
\item In the assignment phase we consider a $K_m$-decomposition of $K_N$ and then describe which components of each $\cC_i$ are assigned to which copies of $K_m$.
\item In the balancing phase we ensure that the mapping from components of each~$\cC_i$ into copies of $K_m$ from $K_N$ will form a $(\xi n)$-balanced packing as promised in Lemma~\ref{l1}.
\end{itemize}
Below we outline the main ideas of these two steps.
We start with the assignment phase first.
The balancing phase will be discussed in Section \ref{outemb}.

\subsection{Outline of the assignment phase}\label{outass}
The purpose of the assignment phase is to produce a ``preliminary packing'' of each $\cC_i$, $i=1,\dots,\bk$ into some $K_m$-factor.
We recall that each component graph $\cC_i$ consists of several components each with at most $s$ vertices and maximum degree at most $\Delta$.
Moreover, in each component $C$ we distinguish the set~$\partial C$ of vertices that are connected to the separator $U_i$.

We define an \emph{isomorphism type} $S$ as a pair $(R,B)$ where $R$ is a graph on at most $s$ labeled vertices and maximum degree at most $\Delta$, and $B$ is a subset of the vertices of $R$.
Let~$\cS=(S_1,\dots,S_\sigma)$ be the enumeration of all isomorphism types $S_j=(R_j, B_j)$, such that
\begin{equation}\label{s2}
\frac{e(R_1)}{v(R_1)}\geq\dots\geq \frac{e(R_\sigma)}{v(R_\sigma)}.
\end{equation}
The definition of $\cS$ yields
\begin{equation}\label{sigma}
\sigma\leq 2^{\binom{s}{2}}\cdot 2^{s}\leq 2^{s^2}.
\end{equation}

For every component $C$ of $\cC_i$ there exists an isomorphism type $S_j=(R_j,B_j)\in\cS$ such that there exists a graph isomorphism $\varphi\colon V(C)\rightarrow V(R_j)$ with the additional property that $\varphi(\partial C)=B_j$.
Therefore, we can describe the structure of a component graph $\cC_i$ as a disjoint union
\[
\cC_i=\bigcup_{S\in\cS}\nu_i(S)\cdot S
\]
where $\nu_i(S)$ denotes the number of components isomorphic to $S$ contained in $\cC_i$.
In the rest of the paper we will simplify the notation and refer to $S$ as a graph.

The assignment procedure makes use of further decomposition layers.
In fact, for each copy of $K_m$ appearing in the resolvable decomposition of $K_N$ we consider a resolvable $K_\ell$-decomposition of such a copy of $K_m$.
Each resolution class consisting of $\frac{m}{\ell}$ disjoint copies of~$K_\ell$ will be reserved for some isomorphism class $S$ and the copies of $S$ coming from various~$\cC_i$ will be then packed into each such $K_\ell$.
Since we consider $K_m$-decomposition of~$K_N$, $K_\ell$-decomposition of $K_m$, and $S$-decomposition of $K_\ell$ for each $S\in\cS$, we will refer to such structure as \emph{three layer decomposition} and motivate its use below.
\subsection{The three layer decomposition}\label{tld}
We begin our discussion with the simpler case when all components in all the component graphs $\cC_i$ are isomorphic to a given graph $S$ and argue why even in this simpler case at least two layers are required.
Then we look at the general case, where the component graphs consist of more different isomorphism types, and explain the use of three layers.

\subsubsection{One layer}
In the case where all components in $\ci$ are isomorphic to a single graph $S$, a straightforward way to pack $\ci$ into $K_N$ would be the following.
Suppose there exists a resolvable $S$-decomposition of $K_N$.
Then, by assigning the components of a graph $\cC_i$ to copies of $S$ from the same $S$-factor, we ensure that the components within each component graph are packed vertex-disjointly.

With this approach, however, we might end up not covering many edges of $K_N$ (and consequently not being able to find a packing of the graphs $\cC_i$).
Let $\cC_1$ and $\cC_2$ be component graphs with strictly more than $N/2$ vertices. 
Once we assign the components of $\cC_1$ to an $S$-factor of $K_N$, we cannot use the other copies of $S$ in the same $S$-factor to accomodate the components of $\cC_2$.
In fact, at least one component of $\cC_2$ would not fit in that $S$-factor and we would have to use a copy of $S$ from another $S$-factor.
We would have to ensure that this copy of $S$ is vertex disjoint from those already used for $\cC_2$ in the previous $S$-factor, and an obvious way to get around this would be to embed all components of $\cC_2$ in a new $S$-factor all together.
However, this would be very wasteful and if many (for example $\Omega(n)$) graphs~$\cC_i$ would be of size strictly larger than $N/2$, then we would not be able to pack all~$\cC_i$ into $K_N$ in such a straightforward way.
We remedy this situation by introducing an additional layer.

\subsubsection{Two layers}
For an appropriately chosen integer $m$, suppose there exist a resolvable $K_m$-decompo\-sition of $K_N$ and a resolvable $S$-decomposition of $K_m$.
Note that, with this additional decomposition layer at hand, we can address the issue raised above more easily.
In fact, we fix a $K_m$-factor of $K_N$ and use sufficiently many $K_m$'s of this $K_m$-factor to host the components of $\cC_1$, all of which are isomorphic to $S$ by our assumption.
The remaining~$K_m$'s of the factor can host the first part of $\cC_2$.
We then ``wrap around'' and reuse the $K_m$'s containing copies of $S$ from $\cC_1$ by selecting a new $S$-factor inside these~$K_m$'s to host the second part of $\cC_2$.
This way the components of $\cC_1$ and $\cC_2$ are packed edge disjointly and the components of $\cC_2$ (resp. $\cC_1$) are in addition vertex disjoint, as required for a packing.
We can continue this process to pack $\cC_3,\cC_4,\dots$ until the fixed $K_m$-factor of~$K_N$ is fully used.
Then we continue with another $K_m$-factor and so on.

This procedure will work if all components of each $\cC_i$ are isomorphic to a single $S$.
Let us note however that in case $\cC_i$ contains components of different isomorphism types two layers may not be sufficient.
This is because we would have to select $S$-factors for different graphs $S$ within $K_m$ and there seems to be no obvious way to achieve this in a two layer decomposition.
Instead we will introduce a third layer, which will give us sufficient flexibility to address this issue.

\subsubsection{Three layers}\label{three}
Here we give an outline and describe how a three layer structure can be used to address the general problem.
The details will follow in section \ref{prfass}.
Consider a resolvable $K_m$-decomposition $\dmn$ of $K_N$, a resolvable $K_\ell$-decomposition $\dlm$ of $K_m$, and resolvable $S$-decompositions $\dsl$ of $K_\ell$ for every $S\in\cS$ (in fact the last assumption will never be used in its full strength, we will use Theorem \ref{R} instead).
We view resolvable decompositions as collections of factors.
We write $\dmn = \{\ddmn_1,\dots,\ddmn_{\frac{N-1}{m-1}}\}$, where~$\ddmn_{\bj}$ is a $K_m$-factor of $K_N$ for $\bj=1,\dots,\frac{N-1}{m-1}$. 

Suppose now we are given graphs $\cC_1,\dots,\cC_{\bk}$, $\cC_i=\bigcup_{S\in\cS}\nu_i(S)\cdot S$.
We will proceed greedily processing the $\cC_i$'s one by one.
In each step we will work with one fixed $K_m$-factor ${\ddmn_{\bj}=\ddmn_\curr}$ of $K_N$ which will be used repeatedly as long as ``sufficiently many'' edges of such factor are available.
For example, $\ddmn_1$ will host $\cC_1,\cC_2,\dots,\cC_a$ for some $a<\bk$, then~$\ddmn_2$ will host $\cC_{a+1},\cC_{a+2},\dots,\cC_b$ for some $a<b<\bk$, and so on.
Once we run out of available edges in factor $\ddmn_\curr$ we will move that factor in the set $\dmn_{\text{used}}\subseteq\dmn$ of factors the edges of which were already assigned to previous $\cC_i$'s and select a new factor $\ddmn_\curr\in\dmn\setminus\dmn_{\text{used}}$ which we will continue to work with.

We outline the assignment within a $K_m$ of the current $K_m$-factor.
For each $K_m\in\ddmn_\curr$ we consider a resolvable decomposition $\dlm=\dlm(K_m)$ of such a $K_m$.
Again some factors in that decomposition might have already been completely used.
Among those which were not completely used yet, we specify $\sigma$ of such ``current'' factors $\ddlm_S$, each ready to be used to embed copies of $S$ in the current particular step.
Since $K_\ell$ admits resolvable $S$-decompositions for every $S\in\cS$, each $\ddlm_S$ corresponds to $\frac{(\ell-1)v(S)}{2e(S)}=t(S)$ $S$-factors of~$K_m$ which we may denote by $\ddslm_1,\dots,\ddslm_{t(S)}$.
At each step, in every $K_m$ we will only use one of such $S$-factors, which we denote by $\ddslm_\curr$.
A set of components of~$\cC_i$ that are going to be assigned to an $S$-factor of a $K_m$ will be referred to as a \emph{chunk}.

With this structure in mind we are able to describe our greedy assignment procedure.
Assume that in the assignment procedure the graphs $\cC_1,\dots,\cC_{i-1}$ were already processed and that $\cC_i=\bigcup_{S\in\cS}\nu_i(S)\cdot S$.
The assignment of $\cC_i$ will consist of the following four steps which we discuss in detail in Section \ref{prfass}.
\begin{enumerate}[label=\rmlabel]
\item For every isomorphism type $S\in\cS$, partition the $\nu_i(S)$ components into as few as possible chunks of size at most $\frac{m}{v(S)}$.
\item For every $S\in\cS$, select $\frac{\nu_i(S)v(S)}{m}$ copies of $K_m$ from the current $K_m$-factor $\ddmn_\curr$ and match each such $K_m$ with a chunk of components isomorphic to $S$.
\item For every $S\in\cS$ and for each chunk of type $S$, assign the components in the chunk to the $S$-factor $\ddslm_\curr$ of $K_m$. The copies of $S$ will cover $m\frac{e(S)}{v(S)}$ edges of $K_m$.
\item Prepare for the assignment of the next component graph.
\end{enumerate}
This procedure leads to a packing of $\ci$ into $K_N$ if we do not run out of $K_m$-factors during the process, and in the proof we shall verify this.
Assuming this for the moment, the procedure above yields a preliminary packing which can be encoded by functions ${\bf=\ffi}$, with $f_i\colon V(\cC_i)\rightarrow V(K_N)$.

\subsection{Outline of the balancing phase}\label{outemb}
In this section we will outline how the preliminary packing $\bf$ obtained in the assignment phase is used to realize a $(\xi n)$-balanced packing of~$\ci$ into $K_N$.
Further detail will be given in Section \ref{prfemb}. 

Note that so far we did not consider the boundary degrees of the vertices of $K_N$ and, in fact, $\bf$ is not guaranteed to be balanced.
However, the layered structure of the assignment will allow us to fix this by using the following degrees of freedom.
Firstly, the $\frac{N}{m}$ $K_m$'s in any of the $\frac{N-1}{m-1}$ $K_m$-factors from $\dmn$ can be permuted independently for each $K_m$-factor.
Since any component graph is assigned to a single $K_m$-factor, the resulting mappings remain injective and the embedding of the $\cC_i$'s stays pairwise edge disjoint.
Secondly, each~$K_m$ can be embedded into $K_N$ in $m!$ possible ways by permuting its vertices.
There are   
\[
\left(\left(\frac{N}{m}\right)!\times(m!)^{\frac{N}{m}}\right)^{\frac{N-1}{m-1}}
\]
such choices in total and each of them leads to a packing of the component graphs $\ci$.

We will pick one of such choices uniformly at random, and show that with positive probability each vertex of $K_N$ is used as a boundary vertex approximately the same number of times.
Since the sum of the boundary degrees is at most $\Delta\delta n^2\leq \xi n^2/2$ (see \eqref{bndr}), this leads to a $(\xi n)$-balanced packing $\bg$ of $\ci$ into $K_N$.

\subsection{Proof of Lemma \ref{l1}}\label{prf}
Given $\xi$ and $\Delta$, set 
\begin{equation}\label{deltal1}
\delta=\frac{\xi}{2\Delta}
\end{equation}
and let $\cG$ be a $(\delta,s)$-separable family, for some $s\in \nat$.
We apply Theorem \ref{R} with 
\begin{equation}\label{etal1}
\eta=\xi/8
\end{equation}
and fix an integer $\ell>s^2$ satisfying that for every $S \in \cS$ there exists an $(S,\eta)$-factorization of $K_\ell$.
Let $m\in\nat$ such that 
\begin{equation}\label{ml1}
m>16\sigma\ell/\xi
\end{equation}
and there exists a resolvable $K_{\ell}$-decomposition of $K_m$ (see Theorem~\ref{RCW}).
Similarly, let
\begin{equation}\label{n0}
n_0>\max\big\{4m^2/\xi,2^{2m}\big\}
\end{equation}
such that for any $n\geq n_0$ satisfying the necessary congruence property there exists a resolvable $K_m$-decomposition of $K_n$.
Having defined $n_0$, we are given a \seq ${\cF=(F_1,\dots,F_{\bk})}$ for some $n\geq n_0$.
We show that there exists a $(\xi n)$-balanced packing of the family of component graphs $\ci$ into $K_N$, for any $N$ with $(1+\frac{\xi}{2})n\leq N\leq(1+\xi)n$ such that $K_N$ admits a $K_m$-decomposition.
Since $n\geq n_0\geq\frac{4m^2}{\xi}$, such $N$ indeed exist.

\subsubsection{The assignment phase}\label{prfass}
Next we elaborate on the outline given in Sections \ref{outass} and \ref{tld}.
First we describe the auxiliary structure we are going to use followed by the actual assignment procedure.

\vspace{1em}

\noindent\textbf{The auxiliary structure.}\quad
For each $S\in\cS$ let $\dsl$ be a fixed $(S,\eta)$-factorization of~$K_{\ell}$ (see Definition \ref{fctr}).
Let $\dlm$ be an arbitrarily chosen resolvable $K_{\ell}$-decomposition of~$K_m$.
Similarly, for the given $N$, denote by $\dmn$ an arbitrarily chosen resolvable $K_m$-decomposition of $K_N$.

At each point of time in the assignment procedure we will work with one $K_m$-factor which we refer to as the \emph{current $K_m$-factor} $\ddmn_\curr \in \dmn$.
Each $K_m$ of the current $K_m$-factor is decomposed into $K_\ell$-factors using $\dlm$.
Moreover, in every $K_m$, for every $S\in\cS$ we pick a $K_\ell$-factor which we denote by $\ddlm_{S}$.
We refer to $\ddlm_S$ as the \emph{current $K_\ell$-factor for $S$}.
We then apply Theorem \ref{R} to all $K_\ell$'s in such a $K_\ell$-factor and obtain $(S,\eta)$-factorizations for every $K_\ell$ in $\ddlm_S$.
Note that we can arbitrarily fix an $S$-matching in each $K_\ell$ of $\ddlm_S$ and obtain an $S$-matching of $K_m$ of size at least 
\begin{equation}\label{dslm}
(1-\eta)\frac{\ell}{v(S)}\frac{m}{\ell}=(1-\eta)\frac{m}{v(S)}.
\end{equation}
This way we set up $t(S)$ edge disjoint $S$-matchings of $K_m$ contained in $\ddlm_S$, for 
\[
(1-\eta)\frac{(\ell-1)v(S)}{2e(S)}\leq t(S)\leq\frac{(\ell-1)v(S)}{2e(S)},
\]
which we denote by $\ddslm_1,\dots,\ddslm_{t(S)}$.
Each of these $S$-matchings cover at least $(1-\eta)\frac{m}{\ell}\binom{\ell}{2}$ edges of the $K_\ell$'s in $\ddlm_S$. 

Every such structure will be used until it is considered \emph{full} according to the following definition.
\begin{dfn}\label{full}
A $K_\ell$-factor $\ddlm_S$ is full when all its $S$-matchings have been used.
A $K_m$ is full when there exists an isomorphism type $S\in\cS$ such that $\ddlm_S$ is full and any other $K_\ell$-factor is either full or reserved to another isomorphism type.
A $K_m$-factor is full when one of its $K_m$'s is full.
\end{dfn}

\vspace{1em}

\noindent\textbf{The assignment procedure.}\quad We now give the details of the four steps outlined in Section \ref{three} for the assignment for the graph $\cC_i=\bigcup_{S\in\cS}\nu_i(S)\cdot S$.

We assume that the graphs $\cC_1,\dots,\cC_{i-1}$ have already been assigned and that the current $K_m$-factor $\ddmn_\curr=\ddmn_{\bj}$ is not full.
\begin{enumerate}[label=\rmlabel]
\item \label{step_i} For each isomorphism type $S\in\cS$ we group the $\nu_i(S)$ copies of $S$ into as few as possible chunks of size at most $(1-\eta)\frac{m}{v(S)}$ (note that this matches the size of an $S$-matching of $K_m$, as given in \eqref{dslm})
The correction factor $(1-\eta)$ here addresses the fact that we deal with $(S,\eta)$-factorizations and not with resolvable $S$-decompositions.
The number $\mu_i(S)$ of chunks required for the $\nu_i(S)$ components of type $S$ is hence given by
\begin{equation}\label{Mi}
\mu_i(S)=\left\lceil\frac{\nu_i(S)\cdot v(S)}{(1-\eta)m}\right\rceil.
\end{equation}

\item \label{step_ii} We order the $K_m$'s in the current $K_m$-factor $\dmn_\curr$ according to the number of edges that have already been assigned to it.
We start with the one in which the least number of edges have been used.
We then assign the $\mu_i(S_1)$ chunks of type~$S_1$ to the first $\mu_i(S_1)$ $K_m$'s in that order and continue in the natural way, that is, the~$\mu_i(S_2)$ chunks of type~$S_2$ are assigned to the next $\mu_i(S_2)$ $K_m$'s, and so on.
Since the members of $\cS$ are ordered non-increasingly according to their densities (see~\eqref{s2}), this way we will ensure that the $K_m$'s in the current $K_m$-factor are used in a balanced way, which is essential to leave only little waste.

\item Once we have determined which chunk goes to which $K_m$, we have to assign the components $S$ of the chunk to their copies in the corresponding $K_m$.
In the chosen~$K_m$ we assign the components of the chunk to $\ddslm_\curr$.
Such a matching exists because we assumed that the current $K_m$-factor $\ddmn_{\bj}$ is not full.
Note that, independently of the precise number of components in the chunk, we use an entire $S$-matching in all the $K_\ell$'s of the current $K_\ell$-factor for $S$ for the assignment of this chunk.
\item \label{step_iv} After we have assigned the components of $\cC_i$ we prepare for the assignment of $\cC_{i+1}$.
In each $K_m$, for every isomorphism type $S$, we check whether the current $K_\ell$-factor for~$S$ is full.
If it is, two cases may arise.
In the first case there exists another $K_\ell$-factor in the $K_m$ that has not been reserved for any $S\in\cS$ yet.
Then, we apply Theorem \ref{R} with $S$ and $\eta$ to all copies of $K_\ell$ in such a $K_\ell$-factor and this factor becomes the current $K_\ell$-factor for $S$, i.e., $\ddlm_S$ in that $K_m$.
In the second case, all $K_\ell$-factors are either full or have been reserved for some $S'\in\cS$ with $S'\neq S$, hence we cannot set up a new $K_\ell$-factor for $S$.
This implies that the $K_m$ and the $K_m$-factor are full (see Definition \ref{full}).
Since we assigned the components of $\cC_i$ to the least used $K_m$'s in the $K_m$-factor, we are ensured that at this point all the~$K_m$'s in $\ddmn_\curr$ are almost completely used.
At this point we add $\ddmn_\curr$ to $\dmn_{\text{used}}$ and set $\ddmn_\curr=\ddmn_{\bj+1}$. 
\end{enumerate}

\vspace{1em}

\noindent\textbf{The assignment phase yields a packing.}\quad
We shall verify that the procedure yields a correct assignment.
For that we have to show that any component graph $\cC_i$ ``fits'' into~$K_N$, and that we do not run out of $K_m$-factors while iterating the four steps for all graphs in~$\ci$.

We first show that every $\cC_i$ fits into one $K_m$-factor.
Recall that in Step~\ref{step_i} the copies isomorphic to some $S\in\cS$ are split into chunks of size at most $(1-\eta)\frac{m}{v(S)}$ and each chunk is assigned to an $S$-matching of $\ddlm_S$.
At this point some vertices may not be used for one of the following two reasons:
\begin{enumerate}[label={\upshape({V\arabic*})}]
\item We always reserve a whole $S$-matching $\ddslm_\curr$ for each chunk, even though some chunks may contain only a few copies of $S$.
In the worst case where only one copy of~$S$ is contained in the chunk we may waste $m-v(S)\leq m$ vertices and in principle this could happen for every isomorphism type $S\in\cS$.
However, since such a ``rounding error'' occurs at most once for each isomorphism type, we may waste at most $\sigma m$ vertices for this reason.
\item We cannot guarantee that the $S$-matchings which we are using are perfect $S$-factors.
However, from Theorem \ref{R} it follows that each matching covers at least 
\[
(1-\eta)\frac{m}{v(S)}v(S)=(1-\eta)m
\] vertices of $K_m$.
Therefore the number of uncovered vertices in the $K_m$-factor due to this imperfection is at most $\eta m\frac{N}{m}=\eta N$.
\end{enumerate}
Hence $\cC_i$ fits into one $K_m$-factor if we ensure that $v(\cC_i)+\sigma m + \eta N\leq N$, which follows from 
\[
v(\cC_i)+\sigma m + \eta N\leq n+\sigma m + \eta N\leq (1+\frac{\xi}{2})n\leq N,
\]
due to \eqref{etal1}, \eqref{ml1}, and \eqref{n0}.

It is left to show that $\frac{N-1}{m-1}$ $K_m$-factors are sufficient to host all the graphs from $\ci$.
For that, we shall bound the number of unused edges in each $K_m$-factor.
At the point when a $K_m$ becomes full, all its $K_\ell$-factors, except for the current $K_\ell$-factors $\ddlm_S$ for each isomorphism type $S\in\cS$, have been used in the assignment.
This leads to the following cases.
\begin{enumerate}[label={\upshape({E\arabic*})}]
\item The current $K_\ell$-factor $\dlm_S$ for a given isomorphism type $S$ may not have been used at all and hence all its $\binom{\ell}{2}\frac{m}{\ell}$ edges are not used in the assignment.
\item Owing to Theorem \ref{R}, in a used $K_\ell$-factor, up to at most $\eta\binom{\ell}{2}\frac{m}{\ell}$ edges are not covered by the $S$-matchings.
\end{enumerate}
Hence the total number of edges that are not used in a full $K_m$ can be bounded by 
\[
\left(\sigma+\eta\frac{m-1}{\ell-1}\right)\binom{\ell}{2}\frac{m}{\ell}.
\]
It is left to establish a similar estimate for the other $K_m$'s in the $K_m$-factor.
Recall that we declared the whole $K_m$-factor to be full as soon as one $K_m$ was full.
Since all components of any $\cC_i\subseteq F_i$ have bounded maximum degree $\Delta$, in each step up to at most $\frac{m\Delta}{2}$ edges are reserved in any $K_m$ of the current $K_m$-factor.
Owing to the balanced selection of the~$K_m$'s within the current $K_m$-factor (see Step~\ref{step_ii}) we have that the number of used edges over all $K_m$'s in $\dmn_\curr$ differs by at most $\frac{m\Delta}{2}$.
Consequently, the number of unused edges in any~$K_m$ at the point when the $K_m$-factor is declared full is at most 
\[
\left(\sigma + \eta\frac{m-1}{\ell-1}\right)\binom{\ell}{2}\frac{m}{\ell} + \frac{m\Delta}{2}.
\]
Using this estimate for all $\binom{N}{2}/\binom{m}{2}$ of the $K_m$ in the $K_m$-decomposition of $K_N$ leads to a total of unused edges of at most
\[
\left(\sigma\frac{\ell-1}{m-1} + \eta + \frac{\Delta}{m-1}\right)\binom{N}{2}<2\eta\binom{N}{2},
\]
where we used~\eqref{etal1}, \eqref{ml1}, and $\Delta<\sigma$.
Furthermore, since by $N\geq(1+\frac{\xi}{2})n$ we have
\[
\binom{n}{2}+2\eta\binom{N}{2}\leq\binom{N}{2},
\]
we have shown that we do not run out of $K_m$-factors and, hence, the assignment procedure yields a preliminary packing of $\ci$.

For the proof of Lemma \ref{l1} we have to show not only that there exists such a packing but also that there is a balanced one.
This will be the focus of the next phase.

\subsubsection{The balancing phase}\label{prfemb}
In the assignment phase we have constructed a preliminary packing $\bf$ of $\ci$ into the $K_m$-factors of $K_N$ as described in Section \ref{outass}.
We now construct a $(\xi n)$-balanced packing $\bh$ by the following random process consisting of two parts.
Firstly, we randomly permute the $\frac{N}{m}$ $K_m$'s in each $K_m$-factor independently and we will denote the resulting packing by $\bg$.
Secondly, for each $K_m$, we pick a random permutation of its vertices.
As we already noted in Section \ref{outemb}, any such permutation yields a packing of $\ci$ into~$K_N$.

It is left to show that with positive probability each vertex $v$ of $K_N$ has boundary degree with respect to $\bh$ bounded by $\xi n$. 
Recall from Definition \ref{xinbp} that the \emph{boundary degree} with respect to $\bf$ of a vertex $v$ is defined by
\[
\dpf(v)=|\{i\in[\bk]\colon f_i^{-1}(v)\in\partial\cC_i\}|.
\]
For a $K_m$ of the $K_m$-decomposition of $K_N$ and a vertex $v$ of $K_m$ we consider the \emph{relative boundary degree}
\[
\dpf(v,K_m)=|\{i\in[\bk]\colon \text{$f_i$ assigns some components of $\cC_i$ to $K_m$ and }f_i^{-1}(v)\in\partial\cC_i\}|
\]
Clearly, $\sum\dpf(v,K_m)=\dpf(v)$, where the sum runs over all $K_m$ from the $K_m$-decomposition of $K_N$ that contain $v$.
For each $K_m$ we define its \emph{label} as the monotone sequence of the relative boundary degrees of its vertices.
Since these labels of the $K_m$'s consist of relative boundary degrees, such a label is invariant under permutations of the vertices
of a $K_m$ and it is invariant under permutations of the $K_m$'s within its $K_m$-factor.
Moreover, the label of a~$K_m$ is determined by the isomorphism types $S\in\cS$ it hosts, because each type $S$ consists of a labelled graph $R$ with a set of boundary vertices $B$.
Since in the assignment phase we assigned an isomorphism type to a whole $K_\ell$-factor, the number of possible labels is bounded by
$|\cS|^{(m-1)/(\ell-1)}=\sigma^{(m-1)/(\ell-1)}<2^m$ (see \eqref{sigma} and the choice of $\ell>s^2$).

For every $K_m$-factor $\dmn_{\bj}$, let $\alpha_{\bj}(A)$ be the number of $K_m$'s with label $A$ in $\dmn_{\bj}$ and define 
\[
\alpha(A)=\sum_{\bj=1}^{\frac{N-1}{m-1}} \alpha_{\bj}(A).
\]
We call a label \emph{common} if $\alpha(A)\geq\frac{\eta}{2^m}\frac{N(N-1)}{m(m-1)}$ and \emph{rare} otherwise.
Note that the total number of $K_m$'s having a rare label is bounded by $\eta\frac{N(N-1)}{m(m-1)}$, therefore
\begin{equation}\label{eq:sumalpha}
\sum_{A\text{ common}}\alpha(A)\geq (1-\eta)\frac{N(N-1)}{m(m-1)}\,.
\end{equation}

We use these labels to show that each vertex in $K_N$ hosts roughly the same amount of boundary vertices.
For that we first prove that an arbitrary vertex is contained in approximately the \emph{expected} number of $K_m$'s of a given common label.
For a vertex $v$ of~$K_N$ and a common label $A$ we denote by $X^{v,A}$ the number of $K_m$'s containing $v$ that have label~$A$.
Note that $X^{v,A}$ is the sum of $\frac{N-1}{m-1}$ indicator variables $X^{v,A}_{\bj}$, where $X^{v,A}_{\bj}=1$ if the $K_m$ from the $K_m$-factor $\dmn_{\bj}$ adjacent to $v$ has label $A$.
The probability that this happens is then given by $\frac{\alpha_{\bj}(A)}{N/m}$.
By applying Chernoff's inequality ((2.9) in \cite{JLR00}) we obtain
\[
\prob\left(|X^{v,A}-\expec X^{v,A}|>\eta\expec X^{v,A}\right)
<2\exp\left(-\frac{\eta^2\expec X^{v,A}}{3}\right)
<2\exp\left(-\frac{\eta^2}{3}\frac{m}{N}\alpha(A)\right).
\]
Consequently, the probability that one of the common labels appears too many or too few times among the $K_m$'s containing some vertex is bounded by 
\[
\sum_{v\in V(K_N)}\sum_{A \text{ common}}2\exp\left(-\frac{\eta^2}{3}\frac{m}{N}\alpha(A)\right)<N 2^{m+1} \exp\left(-\frac{\eta^3(N-1)}{2^m\cdot3(m-1)}\right)<1\,,
\]
where we used that common labels $A$ are defined through $\alpha(A)\geq\frac{\eta}{2^m}\frac{N(N-1)}{m(m-1)}$ in the first inequality.
Therefore, with positive probability, all vertices are \emph{balanced} in the sense that the occurrences of every common label among the $K_m$'s incident to each vertex roughly agree in proportion with the occurrences of that label in the decomposition.

We fix such permutation of the $K_m$'s and the corresponding numbers $X^{v,A}$ for every vertex $v$ and every label $A$.
Let $\bg$ be the corresponding packing of $\ci$ into $K_N$.
As a consequence, we get that for every vertex $v$ the number of $K_m$'s with common labels attached to it satisfies
\begin{align*}
\sum_{A \text{ common}}X^{v,A} &\stackrel{\phantom{\eqref{eq:sumalpha}}}{\geq}\sum_{A \text{ common}}(1-\eta)\frac{m}{N}\alpha(A)=(1-\eta)\frac{m}{N}\sum_{A \text{ common}}\alpha(A)\\
&\stackrel{\eqref{eq:sumalpha}}{\geq}(1-\eta)\frac{m}{N}(1-\eta)\frac{N(N-1)}{m(m-1)}=(1-\eta)^2\frac{N-1}{m-1}\geq(1-2\eta)\frac{N-1}{m-1}.
\end{align*}
We also obtain an upper bound on the number of $K_m$'s with rare labels for every vertex~$v$
\begin{equation}\label{rl}
\sum_{A \text{ rare}}X^{v,A}=\frac{N-1}{m-1}-\sum_{A \text{ common}}X^{v,A}\leq2\eta\,\frac{N-1}{m-1}.
\end{equation}

Next we show that randomly permuting the vertices of each $K_m$ in the $K_m$-decomposition of $K_N$ for the random packing $\bg$ ensures that the boundary degrees in each $K_m$ are evenly distributed.
Let $\dph(v,A)$ be the sum of the boundary degrees of the vertex $v$ within the~$K_m$'s labelled by $A$ and containing $v$.
Clearly,
\[
	\dph(v)=\sum_A\dph(v,A).
\]
We denote by $A(\bj)$ the $\bj$-th element of the degree sequence $A$ and set $\beta(A)=\frac{1}{m}\sum_{\bj=1}^m A(\bj)$ as the average degree in~$A$.
Since $\alpha(A)$ is the number of $K_m$'s with label $A$ in the $K_m$-decomposition of $K_N$ and~$m\beta(A)= \sum_{\bj=1}^m A(\bj)$ is the sum of the relative boundary degrees of the vertices of such a~$K_m$, for later reference we note
\begin{equation}\label{betaalpha}
\sum_A m\beta(A)\alpha(A)=\sum_A \alpha(A)\sum_{\bj=1}^m A(\bj)=\sum_{i=1}^{\bk}|\partial{\cC_i}|.
\end{equation}

For a moment we ignore the $K_m$'s with rare labels, since owing to \eqref{rl} their contribution will be negligible, and consider only those that have a common label.
We first show that for a vertex $v$ of $K_N$ and a common label $A$, $\dph(v,A)$ is in the range $(1\pm \eta)\beta(A)X^{v,A}$ with high probability.
Let $Y^{v,A}_{\bj}$ be the number of $K_m$'s labelled by $A$ in which $v$ gets boundary degree~$A(\bj)$.
By applying Chernoff's inequality we obtain
\[
\prob\left(\left|Y^{v,A}_{\bj}-\frac{X^{v,A}}{m}\right|>\eta \frac{X^{v,A}}{m}\right)<2\exp\left(-\frac{\eta^2}{3}\frac{X^{v,A}}{m}\right)
\]
for every $\bj\in[m]$.
This implies that with probability $1-2m\exp\left(-\frac{\eta^2}{3}\frac{X^{v,A}}{m}\right)$ we have 
\[
\dph(v,A)=\sum_{\bj=1}^{m}A(\bj)Y^{v,A}_{\bj}=\sum_{\bj=1}^{m}A(\bj)(1\pm \eta)\frac{X^{v,A}}{m}=(1\pm \eta)\beta(A)X^{v,A}.
\]

By summing over all common labels, we have that with positive probability there exist permutations for every $K_m$ of the $K_m$-decomposition of $K_N$ for which all vertices have roughly the expected boundary degree.
More precisely, the probability that there exists a misbehaving vertex is bounded by
\begin{align*}
\sum_{v\in V(K_N)}\sum_{A \text{ common}}2m\exp\left(-\frac{\eta^2}{3}\frac{X^{v,A}}{m}\right)
&<N 2^{m+1} m \exp\left(-\frac{\eta^2}{3}(1-\eta)\frac{\alpha(A)}{N}\right)\\
&\leq N 2^{m+1} m \exp\left(-\frac{\eta^3}{2^m\cdot 3}(1-\eta)\frac{N-1}{m(m-1)}\right)<1\,,
\end{align*}
where the first inequality follows from $\bg$ being a packing in which $X^{v,A}$ is close to its expected value for every $v\in V(K_n)$
and the second inequality follows from the definition of common labels.
Therefore, the contribution of the $K_m$'s with common labels for each vertex $v$ is at most
\begin{align*}
\sum_{A \text{ common}}\dph(v,A) &\stackrel{\phantom{\eqref{bndr},\eqref{deltal1}}}{\leq}
\sum_{A \text{ common}}(1+\eta)\beta(A)X^{v,A}\leq\sum_{A \text{ common}}(1+\eta)\left[\beta(A)(1+\eta)\frac{m}{N}\alpha(A)\right]\\
&\stackrel{\phantom{\eqref{bndr},\eqref{deltal1}}}{=}(1+\eta)^2\frac{1}{N}\sum_{A \text{ common}}\left(m \beta(A)\alpha(A)\right)\stackrel{\eqref{betaalpha}}{\leq}(1+\eta)^2\frac{1}{N}\sum_{i=1}^{\bk} |\partial{\cC_i}|\\
&\stackrel{\eqref{bndr},\eqref{deltal1}}{\leq}(1+\eta)^2\frac{1}{N}\frac{\xi}{2} n^2.
\end{align*}

Owing to \eqref{rl}, a vertex can be incident to at most $2\eta\frac{N-1}{m-1}$ $K_m$'s with rare labels.
Since no component of any $\cC_i$ consists of a single isolated vertex (see \ref{sep2} in Definition \ref{deltas}), the largest relative boundary degree of any vertex 
in such a $K_m$ can be at most $m-1$ and we infer that 
\begin{equation}\label{bdfinal}
\dph(v)\leq(1+\eta)^2\frac{1}{N}\frac{\xi}{2} n^2+2\eta\frac{N-1}{m-1}(m-1) < \left(\frac{(1+\eta)^2}{1+\xi/2}\frac{\xi}{2}+2\eta(1+\xi)\right)n 
\stackrel{\eqref{etal1}}{<} \xi n
\end{equation}
for every $v\in V(K_N)$, thus proving Lemma \ref{l1}.
\qed

\section{Packing the separators}\label{L2}
In this section we prove Lemma \ref{l2}.
The Lemma asserts that a balanced packing of~$\ci$ into $K_{(1+\xi)n}$ can be extended to a packing of $\Fi$ in $K_{(1+\eps)}$.
For that we have to show that we can embed the separators $\ui$ in an appropriate way.
Roughly speaking, we will show that a simple greedy strategy will work in here.

\subsection{Proof of Lemma \ref{l2}}
Given $\eps$ and $\Delta$, set 
\[\xi=\frac{\eps}{12\Delta^2}\qquad\text{ and }\qquad\delta=\frac{\eps^2}{72\Delta^2}.
\]
Let $s\in\nat$ and let $\cG$ be a $(\delta,s)$-separable family.
For sufficiently large $n$ let $\cF=(F_1,\dots,F_{\bk})$ be a \seq and suppose that there exists a $(\xi n)$-balanced packing of the component graphs $\ci$ into a clique of order $(1+\xi)n$.
Fix a partition $X\dot{\cup} Y$ of the vertex set of $K_{(1+\eps)n}$, where $|X|=(1+\xi)n$, and denote by $K_X$, $K_Y$, and $K_{X,Y}$ the complete subgraphs induced on $X$ and on $Y$, and the complete bipartite subgraph between $X$ and~$Y$, respectively.
Let $\bh=\hi$ with 
\[
h_i\colon V(\cC_i)\rightarrow X
\]
be a $(\xi n)$-balanced packing of $\ci$ into $K_X$. 
We shall use $K_Y$ to embed $\ui$, and~$K_{X,Y}$ for the necessary connections.
It is easy to see that if the following conditions are satisfied then the resulting map is a packing of $\cF$ into $K_{(1+\eps)n}$:
\begin{enumerate}[label={\upshape({P\arabic*})}]
\item \label{101} for every $i\in [\bk]$, the vertices of $U_i$ are mapped injectively into $Y$;
\item \label{102} each edge in $K_{X,Y}$ is used at most once;
\item \label{103} each edge in $K_Y$ is used at most once.
\end{enumerate}
Note that we will embed $\sum_{i\in[\bk]}|U_i|\leq\delta n^2$ vertices into $Y$, therefore some vertices in $Y$ will be used at least $\frac{\sum_{i\in[n]}|U_i|}{|Y|}\leq\frac{\delta n^2}{|Y|}$ times.
However, we will ensure that each vertex in $Y$ is used at most $3\frac{\delta n^2}{|Y|}$ times.
The packing of $\cF$ into $K_{(1+\eps)n}$ will be expressed by a family of functions $\bhh=\hhi$ with
\[
\hh_i\colon V(F_i)\rightarrow X\dot{\cup}Y
\]
where $\hh_i$ extends $h_i$ from $V(\cC_i)$ to $V(F_i)$.
For a vertex $v\in V(\cC_i)$, we set $\hh_i(v)=h_i(v)\in X$ for any $i\in[n]$.
For the vertices in the separators $\ui$ we will fix their image $\hh_i(v)$ in~$Y$ one by one in a greedy way, starting with vertices of $U_1$.

At each step we embed a vertex $u\in U_i$ into $Y$, assuming that all vertices of $U_j$ with $j<i$ and possibly some (at most $|U_i|-1<\delta n$) vertices of $U_i$ were already embedded.
Let~$N_{\cC_i}(u)$ be the neighbourhood of $u$ in $\cC_i$, and $N_{U_i}(u)$ the neighbourhood of $u$ in $U_i$ both of size at most~$\Delta$.
Suppose so far we made sure that every vertex in $Y$ was used at most~$3\frac{\delta n^2}{|Y|}$ times.
We will embed $u$ in such a way that \ref{101}, \ref{102}, and \ref{103} are obeyed (see~\ref{101x},~\ref{102x}, and~\ref{103x} below), and afterwards each vertex of $Y$ is still used at most~$3\frac{\delta n^2}{|Y|}$ times.
This will show that $\bh$ can be extended to a packing $\bhh$ of $\cF$ and conclude the proof.
Having this in mind we note:
\begin{enumerate}[label={\upshape({P\arabic*$'$})}]
\item \label{101x} The vertices of $U_i$ have to be embedded injectively into $Y$ and, hence, up to at most $|U_i|-1<\delta n$ vertices of $Y$ may not be used for the embedding of $u$.
\item \label{102x} Since every edge in $K_{X,Y}$ can be used at most once, we require $\hh_i(u)\neq\hh_j(u')$ for every vertex $u'\in U_j$ with $\hh_j(N_{\cC_j}(u'))\cap\hh_i(N_{\cC_i}(u))\neq\emptyset$.
Let $x\in\hh_i(N_{\cC_i}(u))$.
Owing to the $(\xi n)$-balancedness of the packing $\hi$, $x$ hosts at most $\xi n$ vertices from $\bigcup_{k\in[\bk]}\partial\cC_k$ and each of them has at most $\Delta$ neighbours in some $U_k$ for $k\in[\bk]$.
Assuming that all of them have already been embedded into $Y$, we obtain at most~$\Delta\xi n$ forbidden vertices for each of the up to at most $\Delta$ neighbours of $u$ in $\cC_i$.
Hence, the total number of forbidden options for $\hh_i(u)$ in $Y$ is at most $\Delta^2\xi n$.
\item \label{103x}
Note that $K_Y$ also hosts the edges contained in the separator $U_i$ and every edge of~$K_Y$ may be used at most once.
Suppose that there exists a vertex $u'$ from $U_j$ with $j<i$ such that $\hh_i(N_{U_i}(u))\cap\hh_j(N_{U_j}(u'))\neq\emptyset$. 
Then $\hh_i(u)$ must avoid $\hh_j(u')$ for any such~$u'$, because at least one edge between this vertex and the image of the neighbours of $u$ is already used.
Since by our assumption every vertex in the set $\hh_i(N_{U_i}(u))$ hosts at most $3\frac{\delta n^2}{|Y|}$ vertices embedded so far, and since $\Delta(F_j)\leq\Delta$, there are at most $\Delta\cdot3\frac{\delta n^2}{|Y|}|N_{U_i}(u)|\leq3\Delta^2\delta n^2/|Y|$ such restrictions.
\end{enumerate}

Since up to now every vertex $y\in Y$ was used at most $3\frac{\delta n^2}{|Y|}$ times for the embedding, by denoting with $Y_u\subseteq Y$ the set of candidates for the embedding of $u$, we obtain 
\[
|Y_u|\geq|Y|-\Big(\delta n + \Delta^2 \xi n + 3 \Delta^2 \frac{\delta n^2}{|Y|}\Big)\geq|Y|-\frac{\eps}{4}n>\frac{|Y|}{2}.
\]
Since we have to embed at most $\sum_{i\in[\bk]}|U_i|\leq\delta n^2$ vertices in total, at any time some vertex $y\in Y_u$ was used at most 
\[
\frac{\delta n^2}{|Y|/2}<3\frac{\delta n^2}{|Y|}-1
\]
times, and this vertex we choose for $\hh_i(u)$.
We have thus shown that at each round we can always pick one vertex in $Y$ such that all the edges needed to connect the vertex we want to embed to all its neighbour are available and it was used before at most $3\frac{\delta n^2}{|Y|}-1$ times.
This completes the proof of the lemma.
\qed

\subsection*{Acknowledgement}
We thank the anonymous referee for her or his constructive remarks.

\end{document}